\documentclass[12pt]{article}
 
\def\w{\dot{w}}
\def\ws{{w_S}}
\def\v{{\rm v}}
\def\sigmad{\dot{\sigma}} 
\def\int{\mathbb{Z}}

\def\Ue{{\cal U}_{\varepsilon}({\mathfrak g})}

\def\OO{{\mathcal O}}

\def\proof{{\bf Proof. }}

\def\pf{\proof}

\def\tw{{\tt w}}
\def\tU{{\tt U}}
\def\tB{{\tt B}}
\usepackage{times}
\usepackage{graphics}
\usepackage{amssymb}
\usepackage{amscd}
\usepackage{amsmath}
\input xy 
\xyoption{all}
\title{A Katsylo theorem for sheets of  spherical conjugacy classes}
\newtheorem{theorem}{Theorem}[section]
\newtheorem{lemma}[theorem]{Lemma}

\newtheorem{proposition}[theorem]{Proposition}

\newtheorem{remark}[theorem]{Remark}

\author{Giovanna Carnovale, Francesco Esposito\\
Dipartimento di Matematica\\
Torre Archimede - via Trieste 63 - 35121 Padova - Italy\\
email: carnoval@math.unipd.it, esposito@math.unipd.it }

\date{}

\begin{document}
\maketitle
\begin{abstract}
We show that, for a sheet or a Lusztig stratum $S$ containing spherical conjugacy classes in a connected 
reductive algebraic group $G$ over an algebraically closed field in good characteristic,
the orbit space $S/G$ is isomorphic to the quotient of an affine subvariety of $G$ modulo
the action of a finite abelian $2$-group. The affine subvariety is a closed subset of a Bruhat 
double coset and the abelian group is a finite subgroup of a  maximal torus of $G$. We  show that 
 sheets of spherical conjugacy classes in a simple group are always smooth and we list which strata containing 
spherical classes are smooth. 
\end{abstract}

\section{Introduction}

In \cite{katsylo}, it is shown that the orbit space of a sheet $S$ of adjoint orbits in a 
complex Lie algebra has the structure of a geometric quotient which is isomorphic to an affine
variety modulo the action of a finite group. The affine variety is the intersection of $S$ with 
the Slodowy slice of a nilpotent element $e$  in $S$, and the finite group is the component group 
of the centralizer of $e$. An algebraic proof of this result was obtained by Im Hof \cite{imhof},
who proved that sheets in complex Lie algebras of classical type are all smooth, by showing that
$S$ is smoothly equivalent to its intersection with the Slodowy slice. 
In addition, Katsylo's quotient is isomorphic to the connected component of Alexeev-Brion invariant Hilbert 
scheme containing the closure of the nilpotent orbit $G\cdot e$ \cite{jansou-ressayre}. Katsylo's theorem has also
been applied to the study of one-dimensional representations of finite $W$-algebras \cite{losev,premet}, 
which is related to the problem of determining the minimal dimensional modules for restricted Lie algebras. 
In this context, it has been shown in \cite{PT} that the space of $1$-dimensional representations of the
finite $W$-algebra associated with a 
nilpotent element $e$ in a classical Lie algebra is isomorphic to an affine space if and only if $e$
lies in a single sheet. The latter condition is equivalent to say that the union of the sheets passing through $e$ is 
a smooth variety. Our goal is to provide an analogue of Katsylo's theorem 
for sheets of conjugacy classes in a reductive algebraic group $G$ over an algebraically closed field of
good characteristic. Since sheets are the irreducible components of the parts in Lusztig's partition 
\cite{lusztig-partition} called strata, the theorem will give an analogue for strata as well. 
As strata should be seen as the group analogue of the union of sheets passing through a nilpotent element $e$, 
used in \cite{PT}, we expect that  their geometry will have relevance in representation theory of 
quantum groups at the roots of unity.

We prove a Katsylo theorem in the case that the sheet (or stratum) in question contains (hence consists of)
spherical conjugacy classes, that is,  classes having a dense orbit for a Borel subgroup $B$ of $G$. 
Strata, and therefore sheets,  of conjugacy classes in a reductive algebraic group do not
necessarily contain unipotent classes, so the analogue of Katsylo's theorem cannot be straightforward. 
A group analogue for Slodowy slices has been introduced in \cite{sevostyanov}. In analogy to Steinberg's 
cross section, these slices depend on a conjugacy class in the Weyl group $W$ of $G$.
The construction of these slices requires a suitable choice of positive roots in the
root system of $G$ which depends on the class of the element in $W$. 
Although the transversality result in \cite{sevostyanov} is stated in characteristic zero, 
the proof holds in arbitrary characteristic. When $w\in W$ acts  without fixed points, a section analogous 
to the one in \cite{sevostyanov} was given in \cite{he-lu}, which contains a 
generalization of Steinberg cross section theorem in this case. 

To our aim, we exploit Sevostyanov's result together with the well-understood behaviour of spherical 
conjugacy classes with respect to the Bruhat decomposition. We replace the Slodowy slice by a suitable subset
${\mathcal S}_w$ of a Bruhat double coset $BwB$, depending on the stratum, such 
that its intersection with each given sheet in the stratum coincides (up to conjugation) 
with the intersection of the sheet with Sevostyanov's slice. Since for spherical conjugacy classes 
the intersection with this double coset is precisely the dense $B$-orbit, 
we show that the intersection of ${\mathcal S}_w$ with 
each conjugacy class is a single orbit for a finite $2$-subgroup of a fixed maximal torus $T$. 

Thanks to Sevostyanov's transversality result, a sheet $S$ of spherical classes is smooth if and only if
$S\cap {\mathcal S}_w$ is so, and similarly for strata. This result is applied in Section \ref{sec:smooth} 
where we obtain the second main result of this paper:  sheets containing a spherical class in simple groups are all smooth. As a consequence, we
classify smooth strata of spherical classes in simple groups. 

\section{Notation}

Unless otherwise stated, $G$ is a connected, reductive algebraic group over an algebraically closed field $k$ of good characteristic, i.e., not bad for any simple component of $[G,G]$. 

Let $T$ be a fixed maximal torus of $G$ and let $\Phi$ be the associated
root system. 
The Weyl group of $G$ will be denoted by $W$.
The centralizer of an element $x\in G$ in a subgroup $H$ of $G$ will be denoted by $H^x$ and its identity component will be denoted by $H^{x\circ}$. 
Let $G$ act regularly on an irreducible
variety $X$.
For $n\geq 0$, we shall denote by $X_{(n)}$ the locally closed subset  $X_{(n)}=\{x\in X~|~ \dim G\cdot x=n\}$. 
For a subset $Y\subset X$, if $m$ is the maximum integer $n$ for which $Y\cap X_{(n)}\neq\emptyset$,  the open subset $X_{(m)}\cap Y$ will be denoted by $Y^{reg}$.  A sheet for the action of $G$ on $X$ is an irreducible component of some $X_{(n)}$. We will investigate the case in which $X=G$ and the action is by conjugation.  Let, for an element $g\in G$, $g=su$ be its Jordan decomposition. It has been shown in \cite{gio-espo} that for any sheet $S$ there is a unique Jordan class $J=J(su)=G\cdot((Z(G^{s\circ})^\circ s)^{reg}u)$ such that $S=\overline{J}^{reg}$. As a set, 
\begin{equation}\label{eq:sheets}S=\bigcup_{z\in Z(G^{s\circ})^\circ}G\cdot(zs {\rm Ind}_{G^{s\circ}}^{G^{zs\circ}}G^{s\circ}\cdot u)\end{equation} where $s$ and $u$ are as above and ${\rm Ind}$ is as in \cite{lusp}.

Sheets of conjugacy classes are the irreducible components of the parts, called strata, of a partition defined in \cite{lusztig-partition} as fibers of a map involving Springer correspondence, \cite{gio-MRL}.

For a Borel subgroup $B\supset T$  with unipotent radical $U$ and a conjugacy class $\OO$ (a sheet $S$, respectively) in $G$, let $w_{\OO}$ ($w_S$, respectively)  be the unique element in $W$ such that $\OO\cap B w_\OO B$ is dense in $\OO$ ($S\cap B w_S B$ is dense in $S$, respectively).  If a sheet $S$ contains a spherical conjugacy class then $w_S=w_\OO$ for every $\OO\subset S$, \cite[Proposition 5.3]{gio-MRL}. In addition, it follows from \cite[Theorem 5.8]{gio-MRL} that $w_S$ is constant along strata containing spherical classes. The element $w_\OO$ is always an involution and it is maximum in its conjugacy class with respect to the Bruhat ordering (\cite{clt,cc}).

The conjugacy classes of $w_\OO$ and $w_S$ in $W$ are independent of the choice of a Borel subgroup containing $T$. Thus, the map $\OO\mapsto w_\OO$ determines a map $\varphi$ from the set of conjugacy classes in $G$ to the set of conjugacy classes of involutions in $W$.

For $w$ an involution in $W$, let $T^w:=\{t\in T~|~w(t)=t\}$, and  $T_w:=\{t\in T~|~w(t)=t^{-1}\}$. Then $T=(T^w)^\circ (T_w)^\circ$ and $S_w:=T_w\cap T^w$ is an elementary abelian $2$-group. For any choice of a Borel subgroup $B$ containing $T$ with unipotent radical $U$, a longest element $w_0\in W$ is determined. 
We set $U^w:=U\cap w^{-1} w_0Uw_0^{-1}w$ and $U_w:=U\cap w^{-1} Uw$.

\section{Spherical classes and Bruhat decomposition}

In this Section we  prove a Katsylo Theorem for sheets containing spherical conjugacy classes. 
We will make use of the following general results.

\begin{lemma}(\cite[Lemma 2.13]{imhof})\label{lem:imhof}Let $A,B,C$ be varieties, let $\eta\colon A\to B$ be a smooth and surjective morphism and let $\theta\colon B\to C$ be a set-theoretic map such that $\theta\eta\colon A\to C$ is a morphism. Then $\theta$ is a morphism. 
\end{lemma}\hfill$\Box$


\begin{lemma}\label{lem:geometric-quotient}Let $X$ be a $G$-variety and assume that there is an affine closed subset $\Sigma\subset X$ with an action of a finite group $\Gamma$ such that the following properties hold:
\begin{enumerate}
\item for every $G$-orbit $\OO$ of $X$ the set $\Sigma\cap\OO$ is a $\Gamma$-orbit;
\item the natural map $\mu\colon G\times\Sigma\to X$ is smooth and surjective. 
\end{enumerate} 
Then, $X/G$ exists and it is isomorphic to $\Sigma/\Gamma$.  
\end{lemma}
\pf Define the map $\psi\colon X\to \Sigma/\Gamma$ set theoretically sending an element $x\in X$ to $G\cdot x\cap \Sigma$. 
We note that $\psi\mu\colon G\times \Sigma\to \Sigma/\Gamma$ is the composition of the natural projection on the second factor followed by the projection to the quotient. 
Therefore Lemma \ref{lem:imhof} applies with $A=G\times\Sigma$, $B=X$, $C=\Sigma/\Gamma$ 
$\eta=\mu$ and $\theta=\psi$, so  $\psi$ is a morphism. 
Since $\mu$ is surjective $\Sigma$ meets all the $G$-orbits in $X$. 
The map $\psi$ is universally submersive because its restriction to $\Sigma$ is the canonical universal quotient map $\pi$ which is universally submersive \cite[Theorem 1.1]{GIT}. 
We now prove that for any open subset $V\subset \Sigma/\Gamma$ we have
$\OO_{\Sigma/\Gamma}(V)\simeq \OO_{X}(\psi^{-1}(V))^G$. It is enough to prove it for $V$ affine. For $f$ a morphism, let  $f^\#$ be the map induced by $f$ on regular functions. Note that $\psi^\#$ embeds $k[V]$ into $k[\psi^{-1}(V)]^G$. Furthermore, if $\iota$ is the inclusion of $\Sigma$ in $X$, then $\iota^\#$ induces an embedding of $k[\psi^{-1}(V)]^G$ into $k[\pi^{-1}(V)]^\Gamma$. Since $\pi^\#=\iota^\#\circ\psi^\#$ induces an isomorphism $k[V]\simeq k[\pi^{-1}(V)]^\Gamma$, also $\iota^\#$ induces an isomorphism $k[\psi^{-1}(V)]^G\simeq k[\pi^{-1}(V)]^\Gamma$, hence $\psi^\#$ induces an isomorphism
$k[V]\simeq k[\psi^{-1}(V)]^G$.  
\hfill$\Box$
\bigskip

Let $\w$ be a representative of $w\in W$ in $N(T)$ and  let $B$ be a Borel subgroup $B$ containing $T$ and with unipotent radical $U$. We define ${\mathcal S}_{w}:=\w T^w U^w$. 

\begin{lemma}\label{lem:affine}Let $S$ be either a sheet or a stratum consisting of spherical conjugacy classes. Then $S\cap {\mathcal S}_{w_{S}}$  is a closed subset of $G$.
\end{lemma}
\pf We show that $S\cap {\mathcal S}_{w_{S}}=\overline{S}\cap {\mathcal S}_{w_{S}}$. 
By the dimension formula in \cite{ccc,gio-mathZ,lu,mauro-cattiva,gio-MRL} we have $S\subset G_{(M)}$ where $M=\ell(w_S)+{\rm rk}(1-w_S)$ and $\overline{S}\setminus S\subset \bigcup_{m<M}G_{(m)}$, \cite[Proposition 5.1]{gio-espo},\cite[Theorem 2.1]{gio-MRL}. If for some class $\OO$ in $G$ and some $w\in W$ there holds $\OO\cap BwB\neq\emptyset$ then $\dim \OO\geq \ell(w)+{\rm rk}(1-w)$, \cite[Theorem 5]{ccc}. Hence, if $\OO\subset \overline{S}\setminus S$ then 
$\OO\cap{\mathcal S}_{w_{S}}\subset \OO\cap B w_S B=\emptyset$. 
\hfill$\Box$
\bigskip

Let us briefly recall the construction of the closed subset $\Sigma_{w}$ of $G$ defined in \cite[p.1890]{sevostyanov}, in the case in which $w\in W$ is an involution. Let $\{v_1,\,\ldots,\,v_r\}$ be a basis of   the $(-1)$-eigenspace of $w$ in the real span  ${\mathfrak h}_{\mathbb R}$  of the co-roots in the Cartan subalgebra ${\mathfrak h}$ and let $\Psi=\Phi\cap {\mathfrak h}_{\mathbb R}^w$. Up to rescaling the $v_i$'s by a positive real scalar, we can construct a set of positive roots $\Phi_+$ satisfying the following rules:
$\Phi_+\cap \Psi$ is defined freely. For $i$ 
maximal satisfying $(\beta,v_i)\neq0$, we have $\beta\in\Phi_+$ if and only if $\beta(v_i)>0$. 

Since $w(v_i)=-v_i$ for every $i$,   there holds $\Phi_+\setminus\Psi=\{\alpha\in\Phi_+~|~w(\alpha)\in -\Phi_+\}$. In other words,  with respect to the constructed choice of positive  roots, $w$ has maximal possible length. In addition, 
$w(\Phi_+\setminus\Psi)=(-\Phi_+)\setminus\Psi$. Let $\tU$ be the subgroup generated by the root subgroups corresponding to roots in $\Phi_+$ and let $\tB:=T\tU$.

Let ${\tw}$ be the unique representative of $w$ such that $\tw x_\alpha(1) \tw^{-1}=x_{w\alpha}(1)$ (\cite[Theorem 5.4.2]{GoGr}), let $L$ be the Levi subgroup of a parabolic subgroup of $G$ containing $T$ and with root system $\Psi$, and let $P^u$ be the unipotent radical of the parabolic subgroup of $G$ containing all root subgroups associated with roots in $(-\Phi_+)\setminus \Psi$. Then $\tw P^u \tw^{-1}=\tU^w$. Sevostyanov's slice in this case is the closed subset 
$$\Sigma_w:=P^u L^{\tw}\tw=P^u \tw L^\tw=\tw \tU^w L^\tw=\tw L^{\tw}\tU^w.$$  

\begin{lemma}\label{lem:t-inv}Let $w\in W$ with $w^2=1$ and let $\w$ be a representative in $N(T)$. 
If $\OO$ is a conjugacy class in $G$ such that $\OO\cap BwB\neq \emptyset$ then $\OO\cap \w (T^w)^\circ U\neq\emptyset$.
\end{lemma}
\pf Clearly $\OO\cap \w TU\neq\emptyset$. Let $x=\w t_wt^w u\in (T_w)^\circ(T^w)^\circ U$. Conjugation by $s\in (T_w)^\circ$ yields $\w s^{-2}t_wt^wu'\in \OO$. Since the square map on $(T_w)^\circ$ is onto, there exists $s\in (T_w)^\circ$ such that $s^{-2}=t_w$, whence the statement. 
\hfill$\Box$

\begin{lemma}\label{lem:intersection} Let $S$ be a stratum or a sheet containing a spherical conjugacy class $\OO$. Let $w\in \varphi(\OO)$.
For $\Phi_+$, $\tB$, $\tU$, $\tw$ and $L$ as in the construction of Sevostyanov's slice and
${\mathcal S}_{w}=\tw T^w\tU^w$, $\Sigma_w=\tw L^{\tw}\tU^w$
 we have 
$$\OO'\cap{\mathcal S}_{w}=\OO'\cap \Sigma_w\neq\emptyset$$ for every class $\OO'\subset S$. 
\end{lemma}
\pf Let $\Delta_+$ be the set of simple roots associated with $\Phi_+$ and let $\Pi:=\Delta_+\cap {\mathfrak h}^w$.
With respect to the given choice of $\Phi_+$ the element $w$ is of maximal length in $\varphi(\OO)$. Hence, it is equal to $w_S$ with respect to the choice of $\tB$. Therefore 
$\OO'\cap \tw \tB\neq\emptyset$ for every $\OO'\subset S$ and, by Lemma \ref{lem:t-inv}, $\OO'\cap \tw T^w\tU\neq\emptyset$. 
By \cite{gio-mathZ}, we have $w=w_0w_\Pi$,  for $w_\Pi$ the longest element in the parabolic subgroup $W_\Pi$ of $W$. Then $L$ is the standard Levi subgroup of the standard parabolic subgroup $P_\Pi$  associated with $\Pi$;  $\tU_w:=\tU_L=\tB\cap L$, and $\tU^{w}=P^u_\Pi$. 

We have ${\mathcal S}_{w}:=\tw T^w\tU^w\subset\tw L^{\tw}\tU^w=\Sigma_w$ so $S\cap{\mathcal S}_{w}\subset S\cap \Sigma_w$. Conversely, 
let $x=\tw lu\in S\cap \tw L^{\tw}\tU^w$. Then $l\in \tB \sigma \tB$ for some $\sigma\in W_\Pi$. Since $\ell(w_0w_\Pi\sigma)=\ell(w_0w_\Pi)+\ell(\sigma)$ for every $\sigma\in W_\Pi$, we have  $x\in \tB w \sigma \tB\cap S$. Maximality of $w_S=w$ among all $\tau$ in $W$ such that $S\cap B\tau B\neq\emptyset$  forces $\sigma=1$.  Hence, $l\in \tB\cap L^{\tw}=T^w\tU_w$ and $S\cap \Sigma_w\subset \tw T^w\tU$. 
Let $\OO'\subset S$ and let  $y=\tw tu\in \tw TU\cap \OO'$. By \cite[Lemmata 4.6, 4.7, 4.8]{gio-mathZ} the only root subgroups occurring in the expression of $u$
are orthogonal to $\Pi$. Hence they lie in $\tU^w$ and 
$\emptyset\neq \OO'\cap \tw T^wU=\OO'\cap \tw T^w \tU^w=\OO'\cap \Sigma_w=\OO'\cap {\mathcal S}_{w}$.\hfill$\Box$

\begin{remark}\label{rem:An}The results contained in \cite{gio-mathZ} and needed in the proof of Lemma \ref{lem:affine} refer to characteristic zero or odd and good. 
However the proofs of Lemmata 4.6, 4.7, 4.8 and Theorems 2.7 and 4.4 therein are still valid for groups of type $A_n$ in characteristic $2$ because also in this case spherical conjugacy 
classes meet only Bruhat cells corresponding to involutions in the Weyl group. This follows  from 
\cite[Theorem 3.4]{mauro-cattiva} for unipotent classes and, in the general case, from results in \cite{mauro-tutte}.
\end{remark}

\bigskip

Let, for an involution $w\in W$
$$\Gamma_w:=\{t\in (T_w)^\circ~|~t^2\in T^w\}=\{t\in T_{w}~|~t^2\in S_w\}.$$ 
 
\begin{theorem}\label{thm:slice}Let $S$ be a stratum or a sheet containing a spherical conjugacy class $\OO$. Let $w\in \varphi(\OO)$. Then  $S/G\simeq (S\cap \Sigma_w)/\Gamma_{w}$.
\end{theorem}
\pf We apply Lemma \ref{lem:geometric-quotient} with $X=S$, $\Sigma=S\cap \Sigma_{w}$ and $\Gamma=\Gamma_w$.

By Lemmata \ref{lem:affine} and \ref{lem:intersection} the set $\Sigma$ is affine and closed in $S$. 
The action of $\Gamma_w$ by conjugation preserves $\Sigma_w$, hence it preserves $\Sigma$. 

Let $\OO$ be a conjugacy class in $S$. we consider $\OO\cap \Sigma=\OO\cap \tw T^w \tU\subset \tB w \tB$. 
Since $\OO$ is spherical and $w=\ws=w_{\OO}$, the set $\emptyset\neq \OO\cap \tB w \tB$ is the dense $\tB$-orbit in $\OO$ by \cite[Theorem 5]{ccc}. 
Therefore, for any $x=\tw t_xu_x$, $y=\tw t_y u_y\in \OO\cap \Sigma$ there is $b=uvs_0s_1\in \tU^w \tU_w (T^w)^\circ (T_w)^\circ$ such that
$uvs_0s_1\tw t_xu_x=\tw t_y u_yuvs_0s_1$. Since $v s_0 s_1\tw\in \tw T\tU$, uniqueness of the Bruhat decomposition $\tU^w w T\tU$ forces $u=1$, 
so $b=vs_0s_1\in \tU_w (T^w)^\circ (T_w)^\circ$. In addition, $(T^w)^\circ\tU_w$ centralizes all elements in $S\cap \tw T^w\tU^w$, \cite[Lemmata 4.6, 4.7, 4.9]{gio-mathZ} so
$y=bxb^{-1}=s_1xs_1^{-1}$, that is, $\tw s_1^{-1}t_xu_x=\tw t_ys_1(s_1^{-1} u_ys_1)$. This implies that $s_1^2=t_xt_y^{-1}\in T^w$ so $s_1\in\Gamma_w$. 

The map $G\times \Sigma_w\to G$ is smooth by \cite[Proposition 2.3]{sevostyanov}. 
The pull-back of this map along the inclusion $S\to G$ is the map $\mu$, and \cite[Theorem III 10.1]{hartshorne} applies.
\hfill$\Box$
\bigskip

\begin{theorem}\label{thm:TU}Let $S$ be a stratum or a sheet containing a spherical conjugacy class $\OO$. Let $B=TU$ be a Borel subgroup of $G$, corresponding to a system of positive roots $\Phi^+$ and a set of simple roots $\Delta$.  Then, for any representative $\w_S$ of $\ws$ we have 
$$S/G\simeq (S\cap {\mathcal S}_\ws)/\Gamma_{w_S}.$$
\end{theorem}
\pf If  $B=\tB$, and $\w_S=\tw$ as in the construction of Sevostyanov's slice, this is Theorem \ref{thm:slice} in force of Lemma \ref{lem:intersection}. 

Let us assume that $\Delta_+\neq\Delta$ and let $\sigma\in W$ such that $\sigma\Delta=\Delta_+$. Then, $w'_S:=\sigma \ws\sigma^{-1}$ is the maximum with respect to the Bruhat ordering determined by $\Delta_+$, i.e., $\tB w'_S\tB\cap S$ is dense in  $S$. Let $\sigmad\in N(T)$ be a representative of $\sigma$. Then $\sigmad T^{\ws}\sigmad^{-1}=T^{w'_S}$ and $\sigmad U^{\ws}\sigmad^{-1}=\tU^{w'_S}$. In addition $\sigmad\w_S \sigmad^{-1}\in \tw'_S (T^{\ws})^\circ (T_{\ws})^\circ$, where $\tw'_S$ is the representative needed for Sevostyanov's construction. Up to multiplying  $\sigmad$  by a suitable element in $(T_{\ws})^\circ$, we can make sure that $\sigmad\w_S \sigmad^{-1}\in \tw'_S (T^{\ws})^\circ$ so
$\sigmad\w_S T^{\ws}U^{\ws}\sigmad^{-1}=\tw'_S (T^{\ws})^\circ  \tU^{w'_S}$. So, conjugation by $\sigmad$ maps $S\cap  {\mathcal S}_\ws$ isomorphically onto
$S\cap \tw'_S T^{w'_S}\tU^{w'_S}=S\cap \Sigma_{w'_S}$ by Lemma \ref{lem:intersection}. A direct verification shows that $\sigmad \Gamma_{\ws}\sigmad^{-1}=\Gamma_{w'_S}$, whence  the statement follows from Theorem \ref{thm:slice}. 
\hfill$\Box$

\section{Smoothness of sheets of spherical classes}\label{sec:smooth}

In this section we detect when $S$ is smooth, for $S$ a sheet or a stratum of spherical classes in a simple group $G$. 
By \cite[Remark 3.4, Proposition 5.1]{gio-espo}  it is enough to consider a representative for each isogeny class of $G$. 
For $G$ of classical type we will use the standard matrix groups. 
For $G$ of exceptional type we shall always consider the simply-connected group.  We will use the classification of spherical conjugacy classes in \cite{gio-pacific}, from which all we have adopted all unexplained notation.

\bigskip

Arguing as in \cite{imhof}, see also \cite[Proposition 3.9 (iii)]{bulois-lie}, we conclude that $S$ is smooth if and only if 
 $S\cap {\mathcal S}_w$ is so. We analyze smoothness of the latter. In order to do so, we recall some information contained in \cite[Lemmata 4.6, 4.7, 4.8]{gio-mathZ} and Remark \ref{rem:An} about this intersection.
 
Let $w_S=w_0w_\Pi$ be the Weyl group element associated with $S$. Then $w_S(\alpha)=\alpha$ for every $\alpha\in \Pi$. Let $\w_S$ be a representative of $w_S$ such that $\w_S x_\alpha(\xi)\w_S^{-1}=x_\alpha(\xi)$ for every $\xi\in k$ and every $\alpha\in \Phi_\Pi$. 
Let $\OO$ be a class in $S$.  If $\w_S t\v\in \OO\cap \w_S TU $ then $\v$ lies in $V_S:=\prod_{\beta\in\Phi^+,\,w_S\beta=-\beta}X_\beta$, where $X_\beta$ is the root subgroup associated with $\beta$.
and $\w_S t$ commutes with $X_\alpha$ for every $\alpha\in{\mathbb Z}\Pi\cap \Phi$. Therefore $S\cap {\mathcal S}_{w_S}\subset \w_S (Z(L_\Pi)\cap T^{w_S})V_S$.

We will make use of the following observation.

\begin{proposition}\label{dixmier}Let $G$ be a simple algebraic group and let $S=\overline{J(su)}^{reg}$ for some $s,u\in G$ be a  sheet of spherical conjugacy classes in $G$. Then either $u=1$ or $S=G\cdot{su}$ and, if $S\neq G\cdot su$, then $S$ contains a semisimple and a unipotent element. \end{proposition}
\pf By the classification of spherical conjugacy classes in \cite{ccc,gio-pacific}, if $rv$ is spherical and $v\neq1$, then $G^{r\circ}$ is semisimple. Therefore, either $G\cdot rv$ is a single sheet, or it lies in $S=\overline{J(s)}^{reg}$, for some semisimple element $s$. In addition, if $S=\overline{J(s)}^{reg}$ is non-trivial, $G^{s\circ}$ is a Levi subgroup, so $S$ contains a unipotent class by \cite[Theorem 5.6(b)]{gio-espo}.\hfill$\Box$ 

\bigskip

We can state the main result of this Section.

\begin{theorem}Let $G$ be a simple algebraic group over $k$.
\begin{enumerate}
\item All sheets of spherical conjugacy classes are smooth. 
\item Let $S$ be a stratum of spherical conjugacy classes. Then $S$ is smooth with the following exceptions:
\begin{itemize}
\item $G$ is of type $B_2$ and $S$ is the stratum containing the unipotent class with partition $(3,1^2)$, or, equivalently, $G$ is of type $C_2$ and $S$ is the stratum containing  the unipotent class with partition $(2^2)$;
\item $G$ is of type $D_{2h+1}$ and $S$ is the stratum containing the unipotent class with partition $(2^{2h},1^2)$.
\end{itemize}
\end{enumerate}
\end{theorem}
\pf By Proposition \ref{dixmier} it is enough to look at sheets containing a semisimple,  element whose connected centralizer os not semisimple.  Their description  follows from \eqref{eq:sheets} and the classification in \cite{ccc,gio-pacific}.  
For all root systems we will  compute the set theoretical intersection of a sheet  $S$ with ${\mathcal S}_{w_S}$ and we will use the following remark. 

\begin{remark}\label{rem:francesco}
The intersection of $S\cap{\mathcal S}_{w_S}$ is reduced. This is proved through the following steps:
\begin{enumerate}
\item $S$ is reduced.
\item The map $G\times({\mathcal S}_{w_S}\cap S)\rightarrow S$ is smooth and surjective, which follows from  
\cite[Proposition 2.3]{sevostyanov} and \cite[III, Theorem 10.1]{hartshorne}.
\item $G\times({\mathcal S}_{w_S}\cap S)$  is reduced; whichfollows from the previous facts using    \cite[\'exp. II, prop.3.1]{SGA}.
\item $G\times({\mathcal S}_{w_S}\cap S)$ is reduced implies that $({\mathcal S}_{w_S}\cap S)$ is reduced.
\end{enumerate}
\end{remark}

\subsection{Type $A_n$}
Let us first consider $H=GL_{n+1}(k)$. In this case sheets and strata coincide, spherical sheets are parametrized by $m=0,\ldots,\left[\frac{n+1}{2}\right]$ and they are as follows:
$S_m=Z(H)^\circ \OO_m\cup \bigcup_{\lambda,\mu\in k^*;\;\lambda\neq\mu}\OO_m(\lambda,\mu)$, where $\OO_k$ is the unipotent class corresponding to the partition $(2^m,1^{n+1-2m})$ and $\OO_m(\lambda,\mu)$ is the semisimple class with eigenvalues $\lambda$ with multiplicity $n+1-2m$ and $\mu$ with multiplicity $2m$, the case $m=0$ being trivial. The Weyl group element associated to $S_m$ is $w_0w_\Pi$ where $\Pi=\{\alpha_{m+1},\alpha_{m+2},\,\ldots,\alpha_{n+1-m}\}$. We choose
$\w_S=\left(\begin{smallmatrix}
0&0&J_m\\
0&I_{n+1-2m}&0\\
-J_m&0&0
\end{smallmatrix}\right)$ where $J_m$ is the $m\times m$ matrix with $1$ on the antidiagonal and $0$ elsewhere. 
Then 
$$\begin{array}{l}
\w_S (Z(L_\Pi)\cap T^{w_S})V_S\\
=\left\{\left(\begin{smallmatrix}
&&&&&&&&a_1\\
&&&&&&&a_2\\
&&&&&&\dots\\
&&&&&a_m\\
&&&&bI_{n+1-2m}\\
&&&-a_m&&a_m\zeta_m\\
&&\dots&&&&\dots\\
&-a_2&&&&&&a_2\zeta_2\\
-a_1&&&&&&&&a_1\zeta_1
\end{smallmatrix}\right), a_i,b\in k^*;\;\zeta_i\in k\right\}.\end{array}$$
A matrix in $\w_S (Z(L_\Pi)\cap T^{w_S})V_S$ lies in $S_m$ if either all its eigenvalues are equal or else it has  two eigenvalues and it is semisimple. This happens if and only if 
$$\begin{array}{lr}Tr\left(\begin{smallmatrix}0&a_i\\
-a_i&-a_i\zeta_i\end{smallmatrix}\right)=Tr\left(\begin{smallmatrix}0&a_j\\
-a_j&-a_j\zeta_j\end{smallmatrix}\right),& \forall i,j\\
\det\left(\begin{smallmatrix}0&a_i\\
-a_i&-a_i\zeta_i\end{smallmatrix}\right)=\det\left(\begin{smallmatrix}0&a_j\\
-a_j&-a_j\zeta_j\end{smallmatrix}\right);&\forall i,j\\
 b^2+a_1\zeta_1b+a_1^2=0,
\end{array}$$
that is, if and only if there exist $\epsilon_2,\ldots,\epsilon_m\in\{0,1\}$ such that
$a_i=\epsilon_ia_1$, $\zeta_i=\epsilon_i\zeta_1$ and $\zeta_i=-a_1^2b^{-1}-b^2a_1^{-1}$.
The set theoretic intersection $S_m\cap {\mathcal S}_{\w_S} $ is then a union of $(m-1)$ disjoint irreducible components each isomorphic to the image of the morphism
$$\begin{array}{rl}
f\colon k^*\times k^*&\to k^*\times k\cr
(a,b)&\mapsto (a,b,a^2b^{-1}+b).
\end{array}$$
Being a graph, this intersection is smooth for every field $k$ and every $m$. By Remark \ref{rem:francesco}, this intersection coincides with the scheme theoretic one. 

Let us now consider $G=SL_n(k)$. Set theoretically, every sheet of spherical classes $S$ is contained in the intersection of some $S_m$ with $G$.
If ${\rm char}(k)=p$ does not divide $n+1$, then $S\cap{\mathcal S}_{w_S}$ is contained in the image through $f$ of the disjoint smooth curves ${\mathcal C}_{\pm1}$ of equation $a^{2m}b^{n+1-2m}-\pm1=0$. By Remark \ref{rem:francesco}, this set theoretic inclusion is scheme theoretic hence  $S\cap {\mathcal S}_{w_S}$  is smooth.  

Let us now assume that $p|n+1$. Then, for every $m$ coprime with $p$, the argument above applies. If, instead, $p|m$ then 
we still have the set-theoretical inclusion $S\cap {\mathcal S}_{w_S}\subset f({\mathcal C}_{1}\cup {\mathcal C}_{-1})$ but
the curves ${\mathcal C}_{\pm1}$ are not reduced. The reduced scheme of $f({\mathcal C}_{1}\cup {\mathcal C}_{-1})$ is smooth and the above argument applies.

\subsection{Type $B_n$}

Let $G=SO_{2n+1}(k)$ with $n\geq 2$. The non-trivial sheets of spherical conjugacy classes  are given by $S$ and $S'$, with
$$S=\left(\cup_{\lambda\neq0,\pm1}\OO_{\lambda}\right) \cup \OO_{(3,2)}\cup G\cdot \rho_nu$$ where $\OO_{\lambda}$ is the semisimple class with eigenvalues $1,\lambda,\lambda^{-1}$ with multiplicity $1,n,n$ 
respectively;
$\OO_{(3,2)}$ is the unipotent conjugacy class corresponding to the partition $(3,2^{n-1})$, for $n$ odd and $(3,2^{n-2},1^2)$
for $n$ even; the element $\rho_n$ is the isolated diagonal matrix
${\rm diag}(1,-I_{2n})$ and $u$ is a representative of any unipotent conjugacy class in 
$G^{\rho_n\circ}\cong SO_{2n}(k)$ associated with the
partition $(2^{n})$ when $n$ is even, and $(2^{n-1},1^2)$ when $n$ is odd; and 
$$S'=\left(\cup_{\lambda\neq0,1}\OO_{\lambda,1}\right) \cup \OO_{(3,1^{2n-2})},$$
where $\OO_{\lambda,1 }$ is the class of a semisimple matrix with eigenvalues 
$1,\lambda,\lambda^{-1}$ with multiplicity $2n-1,1,1$, respectively and 
$\OO_{(3,1^{2n-2})}$ is the unipotent conjugacy class with associated partition $(3,1^{2n-2})$. 

\bigskip

We have $S\cap S'=\emptyset$ unless $n=2$, so the stratum containing $S$ is not smooth for $n=2$ whereas for $n\geq 3$ the strata containing $S$ and $S'$ are smooth if and only if $S$ and $S'$ are so.

\bigskip

Let us analyze $S$. Here, $w_{S}=w_0$. If we choose $\w_{S}=\left(\begin{smallmatrix}
(-1)^n&0&0\\
0&0&I_n\\
0&I_{n}&0\end{smallmatrix}\right)$, then 
${\mathcal S}_{w_S}=\w_{S} T^{w_{S}}U$
consists of matrices of the form
$$X=X(E,M,Q,v)=\left(\begin{smallmatrix}
(-1)^n&0&(-1)^n\,^t\!v\\
0&0&E\,^t\!Q^{-1}\\
-EQv&EQ&EQM\end{smallmatrix}\right)$$ where $E\in \{\pm1\}^n$,
$v=^t\!(v_1,\,\ldots,\,v_n)\in k^n$, $Q$ is a unipotent upper triangular matrix in $GL_n(k)$, and $M=(-1/2)v\,^t\!v+A$, where $A$ is skew-symmetric.

Now, if $X$ lies in $S$ then there exists $\lambda\in k^*$ such that ${\rm rk}(X-\lambda I)\leq n+1$. 

Assume first that $X=X(E,M,Q,v)$ satisfies ${\rm rk}(X-\lambda I)\leq n+1$ for some $\lambda\neq (-1)^n$.
Then we have
\begin{equation}\label{eq:globale}M-\lambda Q^{-1}E+\lambda^{-1}E^t\!Q^{-1}+\frac{(-1)^n}{(-1)^n-\lambda}v\,^t\!v=0.\end{equation}
Let $\varphi_{\lambda,n}:=\frac{(-1)^n+\lambda}{2((-1)^n-\lambda)}$. Taking symmetric and skew-symmetric parts in \eqref{eq:globale} we obtain the following equations:
\begin{equation}\label{eq:simmetrica}
\varphi_{\lambda,n} v\,^t\!v=(1/2)(\lambda-\lambda^{-1})(Q^{-1}E+E^{^t}\!Q^{-1})
\end{equation}
and
\begin{equation} 
\label{eq:antisimmetrica}
A=(1/2)(\lambda+\lambda^{-1})(Q^{-1}E-E^{^t}\!Q^{-1}).
\end{equation} 
The diagonal terms in \eqref{eq:globale} give
\begin{equation}\label{eq:v-lambda}\varphi_{\lambda,n}e_i v_i^2=(\lambda-\lambda^{-1}).\end{equation}
Hence, if in addition $\lambda\neq(-1)^{n+1}$, then $e_iv_i^2=e_1v_1^2$ for every $i$.
We fix, for $i=1,\,\ldots, n$,  elements $\zeta_i\in k$ such that $\zeta_i^2=e_i$ and we set $a_\lambda:=\zeta_1 v_1$, so $\varphi_{\lambda,n}a_\lambda^2=(\lambda-\lambda^{-1})$. Therefore, for every $j\geq 1$ there is $\eta_j=\pm 1$, with $\eta_1=1$, such that $\zeta_jv_j=\eta_j a_\lambda$. Thus, we have
\begin{equation}\label{eq:alambda}\lambda^2-(2(-1)^n-a_\lambda^2/2)\lambda+1=0
\end{equation}
which gives 
\begin{equation}\label{eq:lambda}(\lambda+\lambda^{-1})=(2(-1)^n-a_\lambda^2/2).
\end{equation}

Making use of \eqref{eq:v-lambda} and \eqref{eq:lambda}, for $2\leq i<j\leq n$, the $(i,j)$ 
entries of \eqref{eq:simmetrica}  and \eqref{eq:antisimmetrica}  give
\begin{equation}\label{eq:QA}(Q^{-1})_{ij}=2\zeta_i^{-1}\zeta_j\eta_i\eta_j,\quad a_{ij}=\left(2(-1)^n-a_\lambda^2/2\right)\zeta_i^{-1}\zeta_j^{-1}\eta_i\eta_j.\end{equation} 

So, for $\lambda\neq\pm1$ and for every choice of $\eta_i$, $\zeta_i$, the matrix $Q$ is completely determined, the vector $v$ depends linearly on $a_\lambda$ and $M$ depends on $a_\lambda^2$, giving a dense subset of a line.  Conversely, if $\lambda\neq\pm1$, the condition ${\rm rk}(X-\lambda I)\leq n+1$ 
also implies that $X$ is semisimple, and it ensures $X\in S$.

Let us now assume $\lambda=(-1)^{n+1}$. Then, \eqref{eq:antisimmetrica} gives
$a_{ij}=(-1)^{n+1}(Q^{-1})_{ij}e_j$ for every $i<j$. 
$X$ lies in $S$ only if ${\rm rk}(X-(-1)^{n+1})^2=1$. Looking at the $(2,2)$-block in this matrix we get 
${\rm rk}(^t\!Q^{-1}E+EQ^{-1})\leq1$, which yields
$$(Q^{-1})_{ij}=2\zeta_i^{-1}\zeta_j\eta_i\eta_j,\text{ and }a_{ij}=2(-1)^{n+1}\zeta_i^{-1}\zeta_j^{-1}\eta_i\eta_j.$$
Let $N=(X-(-1)^{n+1})^2{\rm diag}(1,Q^{-1}E,I_n)$. 
Every row in $N$ must be  a multiple of the first one, which is nonzero. 
Thus, every row in the block $(2,2)$ must be a multiple of the $(1,2)$-block. This gives 
$e_iv_i^2=e_1v_1^2\neq0$ for every $i\geq0$. We set $a=\zeta_1 v_1$, so  for every $i\geq2$ we have $v_i=a\zeta_i^{-1}\eta_i$. A direct computation shows that 
$^t\!v EQv=0$ and $-E\,^t\!Q^{-1}EQv=v$. The condition that the principal minor of size $2$ must be $0$ gives $a^2=8(-1)^n$, i.e.,  $a$ satisfies condition \eqref{eq:alambda} so 
\eqref{eq:QA} is verified also in this case. Thus, for the Jordan class $J$ of ${\rm diag}(1,\lambda I_n,\lambda^{-1} I_n)$, 
 the set-theoretical intersection $J\cap {\mathcal S}_{w_0}$ is a disjoint union of $2^{2n-1}$ copies of $k^*$, given by the values of $a_\lambda$, 
 one for each choice of each $e_j$'s and of the $\eta_j$'s. 

Let us now consider $\lambda=(-1)^n$. Then $X$ must satisfy  
the condition ${\rm rk}(X-(-1)^{n})=n$, which forces $a_\lambda=0$. 
Moreover, $X$ lies in $S$ only if ${\rm rk}(X-(-1)^{n})=1$, which,
combined with \eqref{eq:antisimmetrica} gives the condition 
${\rm rk}(^t\!Q^{-1}E+EQ^{-1})\leq1$, which yields \eqref{eq:QA} with $a_\lambda=0$. 
By the dimensional argument used in the proof of Lemma \ref{lem:affine} for the $G$-conjugacy class of  $X$ shows that $X\in S$.
Hence the set theoretical intersection $S\cap {\mathcal S}_{w_S}$ is a disjoint union of $2^{2n-1}$ copies of $k$.  
By Remark \ref{rem:francesco}, this is also the scheme theoretic intersection so $S$ is smooth.

\bigskip

Let us now consider $S'$. In this case $w_{S'}=s_{\beta}s_1=w_0w_\Pi$ for $\Pi=\{\alpha_3,\ldots,\alpha_n\}$ and $\beta=\varepsilon_1+\varepsilon_2$ the highest long root. 
We choose $\w_{S'}=\left(\begin{smallmatrix}
1&0&0&0&0\\
0&0&0&I_2&0\\
0&0&I_{n-2}&0&0\\
0&I_2&0&0&0\\
0&0&0&0&I_{n-2}\end{smallmatrix}\right)$ so 
$$\begin{array}{l}
\w_{S'} (Z(L_{\Pi})\cap T^{w_{S'}})V_{S'}\\
=\bigcup_{\epsilon,\eta=\pm1}\left\{\left(\begin{smallmatrix}
1&0&0&0&a&b&0\\
0&0&0&0&\epsilon&0&0\\
0&0&1&0&-\eta x&\eta&0\\
0&0&0&cI_{n-2}&0&0&0\\
-\epsilon(a+bx)&\epsilon&\epsilon x&0&\epsilon l&\epsilon m&0\\
-\eta b&0&\eta&0&\eta m'&-\frac{1}{2}\eta b^2&0\\ 
0&0&0&0&0&0&c^{-1}I_{n-2}
\end{smallmatrix}\right)
\begin{array}{l}
a,\,b,x,l,m,\in k; c\in k^*,\\
 m'=-m-ab--\frac{1}{2}\eta b^2\;\end{array}
\right\}.\end{array}$$
Then an element $X\in\w_{S_1} (Z(L_{\Pi_1})\cap T^{w_{S_1}})V_{S_1}$ lies in $S'$ if and only if $rk(X-I)=2$: this is clear if the eigenvalues different from $1$ are distinct. If the eigenvalues different from $1$ are equal to $-1$, then this follows from the fact that the unipotent part must lie in the connected centralizer of the semisimple part. If the eigenvalues are all equal to $1$, then there are only two unipotent classes for which ${\rm rk}(X-I)=2$, namely the one associated with $(2^2, 1^{2n-3})$ and  $\OO_{(3,1^{2n-3})}$. By dimensional reason, the former does not intersect ${\mathcal S}_{w_s'}$.

Assume $rk(X-I)=2$. For such an $X$ we have
$$c=1,\quad a=b=0, \quad l=\eta x^2, m=-\eta x.$$
By Remark \ref{rem:francesco} the variety $S'\cap {\mathcal S}_{w_{S'}}$ is isomorphic to a disjoint union of $4$ affine lines, one for each value of $\eta$ and $\epsilon$.

\subsection{Type $C_n$}

Let us consider $G=Sp_{2n}(k)$ with $n\geq3$. 
There are, up to a central element, two non-trivial sheets of spherical classes, $\pm S_1$ and $S_2$ where
$$S_1=\left(\cup_{\lambda\neq0,\pm1}\OO_{(\lambda,1)}\right)\cup \OO_{(2^2,1^{2n-4})}\cup G\cdot \OO_{\sigma_1 x_{\beta}(1)}$$
where $\OO_{(\lambda,1)}$ is the semisimple class with eigenvalues  $\lambda$, $\lambda^{-1}$ and $1$,
with multiplicity $1,1,2n-2$ respectively, $\OO_{(2^2,1^{2n-4})}$ is the unipotent conjugacy class corresponding 
to the partition $(2^2,1^{2n-4})$,  the element $\sigma_1$ is the isolated 
diagonal matrix ${\rm diag}(-1,I_{n-1}, -1, I_{n-1})$ and $\beta=\varepsilon_1+\varepsilon_2$ is the highest root;
and
$$S_2=\left(\cup_{\lambda\neq0,\pm1}\OO_{\lambda}\right)\cup \pm \OO_{(2^n)}$$
where $\OO_{\lambda}$ is the semisimple class with eigenvalues $\lambda^{\pm1}$
and $\OO_{(2^n)}$ is the unipotent conjugacy class corresponding to the partition $(2^n)$.

\bigskip

Since $n>2$, we always have $S_1\cap S_2=S_2\cap (-S_1)=\emptyset$, hence the strata containing these sheets
are smooth if and only if the sheets are so. 

\bigskip

The Weyl group element corresponding to $S_1$ is $w_{S_1}=w_0w_{\Pi_1}$, for $\Pi_1=\{\alpha_3,\,\ldots,\alpha_n\}$, so $w_{S_1}=s_{\alpha_1}s_\beta$.
We choose $\w_{S_1}=\left(\begin{smallmatrix}
&&I_2\\
&I_{n-2}\\
-I_2\\
&&&I_{n-2}\end{smallmatrix}\right)$ so 
$$\begin{array}{l}
\w_{S_1} (Z(L_{\Pi_1})\cap T^{w_{S_1}})V_{S_1}\\
=\bigcup_{\epsilon,\eta=\pm1}\left\{\left(\begin{smallmatrix}
&&&\epsilon&0\\
&&&-\xi\eta&\eta\\
&&bI_{n-2}\\
-\epsilon&-\epsilon\xi&&-\epsilon(x+\xi y)&-\epsilon(y+\xi z)\\
0&-\eta&&-\eta y&-\eta z\\
&&&&&b^{-1}I_{n-2}
\end{smallmatrix}\right),b\in k^*;x,y,z,\xi\in k\right\}.\\
\end{array}$$
Then $X\in\w_{S_1} (Z(L_{\Pi_1})\cap T^{w_{S_1}})V_{S_1}$ lies in $S_1$ if and only if $rk(X-I)=2$, which holds if and only if
$$b=1,\quad x=-2\epsilon,\quad, y=\eta\xi,\quad z=-2\eta$$
By Remark \ref{rem:francesco} the variety $S_1\cap {\mathcal S}_{w_{S_1}}$ is isomorphic to a disjoint union of affine lines.

\bigskip

The Weyl group element corresponding to $S_2$ is $w_0$ and we choose the representative $\w_0=\left(\begin{smallmatrix}0&I_n\\
-I_n&0\end{smallmatrix}\right)$.
Then the matrices in $\w_0 T^{w_0}V_{S_2}$ are all matrices of the form
$x(E,V,X)=\left(\begin{smallmatrix}0&E^tV^{-1}\\
-EV&-EV X\end{smallmatrix}\right)$, where $E={\rm diag}(\epsilon_1,\,\ldots, \epsilon_n)$, $\epsilon_i=\pm1$ for every $i$, $V$ is an upper triangular unipotent matrix, and $X$ is a symmetric matrix. 

If $x(E,V,X)$ lies in $S_2$, then there is a $\lambda\in k^*$ for which ${\rm rk}(x(E,V,X)-\lambda I)=n$. 
This forces $\lambda X+\lambda^2(V^{-1}E)+^t\!(V^{-1}E)=0$. If $\lambda^2\neq 1$ this can happen only if $V=I$ and $X=-(\lambda+\lambda^{-1})E$ and 
$x(E,V,-(\lambda+\lambda^{-1})E)\in\OO_{\lambda}$. If instead $\lambda^2=1$, then $x(E,V,X)$ lies in $S_2$ only if $(x(E,V,X)-\lambda I)^2=0$. A direct computation shows that this is possible only if $V=I$, so $X=-2\lambda^{-1}E=-(\lambda+\lambda^{-1})E$. If this is the case, $x(E,I,-2\lambda^{-1}E)\in \lambda\OO_{(2^n)}$. Therefore, set theoretically, $S_2\cap {\mathcal S}_{w_0}$ is a disjoint union of affine lines with coordinate $\lambda+\lambda^{-1}$. We apply Remark \ref{rem:francesco}. 
 
\subsection{Type $D_n$}
Let $G=SO_{2n}(k)$, for $n\geq 4$. It is convenient to separate the cases of $n$ even and odd.

\subsubsection{$D_n$ for $n$ even}

Let $n=2h$. Let $\theta$ be a non-trivial graph automorphism of $G$. The only non-trivial sheets of spherical conjugacy classes are given by $S$, $\theta(S)$ and $S'$ as follows:
$$S=\cup_{\lambda\neq0,\pm1}\left(G\cdot\begin{pmatrix}
\lambda I_n\\
&\lambda^{-1} I_{n}\\
\end{pmatrix}\right)\cup \pm \OO_{(2^n)}$$ and
$$\theta(S)=\cup_{\lambda\neq0,\pm1}\left(G\cdot\begin{pmatrix}
\lambda I_{n-1}\\
&\lambda^{-1} I_{n}\\
&&\lambda
\end{pmatrix}\right)\cup \pm \OO'_{(2^n)}$$
where $\OO_n$ and $\OO_n'$ are the two distinct unipotent conjugacy class corresponding to the partition $(2^n)$ and
$$S'=\left(\cup_{\lambda\neq0,1}\OO_{\lambda,1}\right) \cup \OO_{(3,1^{2n-3})}$$
where $\OO_{\lambda,1 }$ is the class of a semisimple matrix with eigenvalues $1,\lambda,\lambda^{-1}$ with multiplicity $2n-2,1,1$, respectively and 
$\OO_{(3,1^{2n-3})}$ is the unipotent conjugacy class with associated partition $(3,1^{2n-3})$. 

The intersection of any pair of distinct sheets is trivial, and the stratum  is smooth if and only if the sheets it contains are so. 

It is enough to deal with $S$ and $S'$. 

For the sheet $S$ we have $w_S=s_{\varepsilon_1+\varepsilon_2}s_{\varepsilon_3+\varepsilon_4}\cdots s_{\varepsilon_{n-1}+\varepsilon_n}=w_0w_\Pi$ for  $\Pi=\{\alpha_1,\alpha_3,\ldots,\alpha_{n-1}\}$. In this case $\theta(w_S)\neq w_S$. We choose the representative
$\w_S=\left(\begin{smallmatrix}
0&L\\
L&0
\end{smallmatrix}\right)$ where $L={\rm diag}(J,J,\ldots,J)$ and $J=\left(\begin{smallmatrix}
0&1\\
-1&0
\end{smallmatrix}\right)$. 

Then, for  $i=1,\,\ldots, h$ and $\epsilon_i=\pm1$,  $\w_{S} (Z(L_{\Pi})\cap T^{w_{S}})V_{S}$ is the disjoint union of the sets of matrices of the form
$x(E,D)=\left(\begin{smallmatrix}
0&E\\
E& D\end{smallmatrix}\right)$ with $x_i\in k$ for $i=1,\ldots, h$ and 
$$E={\rm diag}(E_1,\ldots,E_h),\quad E_i=\left(\begin{smallmatrix}
0&\epsilon_i\\
-\epsilon_i& 0
\end{smallmatrix}\right),\quad D={\rm diag}(-\epsilon_1 x_1 I_2,\ldots,-\epsilon_h x_h I_2).$$
Then $x(E,D)$ lies in $S$ only if there exists $\lambda\in k^*$ such that ${\rm rk}(x(E,D)-\lambda I)=n$. This is possible only if $D=(\lambda+\lambda^{-1})I$. Conversely, if this is the case, a direct verification shows that $x(E,D)$ is either semisimple with eigenvalues $\lambda^{\pm1}$ or unipotent up to a sign. In addition, as $w_{\theta(S)}\neq w_S$ have the same length, and $w_{\theta(S)}$ is maximal among the elements $\tau\in W$ such that $\theta(S)\cap B\tau B\neq\emptyset$, we see that $\OO\cap B w_S B=\emptyset$ for every $\OO\subset \theta(S)$. Thus, any $x(E,D)$ satisfying  ${\rm rk}(x(E,D)-\lambda I)=n$ lies in $S$. 
So, set theoretically, $S\cap {\mathcal S}_{w_S}$ is a disjoint union of $2^h$ affine lines. We conclude as in the previous cases.

For the sheet $S'$ we have $w_{S'}=s_{\beta}s_1=w_0w_\Pi$ for $\Pi=\{\alpha_3,\ldots,\alpha_n\}$ and $\beta=\varepsilon_1+\varepsilon_2$ the highest root. 
We choose $\w_{S'}=\left(\begin{smallmatrix}
&&I_2\\
&I_{n-2}\\
I_2\\
&&&I_{n-2}\end{smallmatrix}\right)$ so 
$$\begin{array}{l}
\w_{S'} (Z(L_{\Pi})\cap T^{w_{S'}})V_{S'}\\
=\bigcup_{\epsilon,\eta=\pm1}\left\{\left(\begin{smallmatrix}
0&0&0&\epsilon&0&0\\
0&1&0&-\eta x&\eta&0\\
0&0&cI_{n-2}&0&0&0\\
\epsilon&\epsilon x&0&\epsilon l&\epsilon m&0\\
0&\eta&0&\eta-m&0&0\\ 
0&0&0&0&0&c^{-1}I_{n-2}
\end{smallmatrix}\right),
x,l,m,\in k; c\in k^*,
\right\}.\end{array}$$

If $X\in\w_{S'} (Z(L_{\Pi_1})\cap T^{w_{S'}})V_{S'}$ lies in $S'$ then $rk(X-I)=2$. All elements satisfying this condition lie in $S'$. Indeed,  the centralizer of the representatives of the classes in $S'$ in $O_{2n}(k)$ is not contained in $SO_{2n}(k)$ so elements that are $GL_{2n}(k)$-conjugate, are also $SO_{2n}(k)$-conjugate. Therefore, 
the argument used for the sheet $S'$ in type $B_n$ applies. For such an $X$ we have
$$c=1,\quad l=\eta x^2, m=-\eta x.$$
Hence the variety $S'\cap {\mathcal S}_{w_{S'}}$ is isomorphic to a disjoint union of $4$ affine lines.

\subsubsection{$D_n$ for $n$ odd}
Let $n=2h+1$. The only non-trivial sheets of spherical conjugacy classes are $R$, $\theta(R)$ and $S'$ as follows: 
$$R=\cup_{\lambda\neq0,\pm1}\left(G\cdot\begin{pmatrix}
\lambda I_n\\
&\lambda^{-1} I_{n}\\
\end{pmatrix}\right)\cup \pm \OO_{(2^{n-1},1^2)};$$
where $\OO_{(2^{n-1},1^2)}$ is the unique unipotent conjugacy class corresponding to the partition $(2^{n-1},1^2)$; and 
$S'$  is the same as for $n$ even and can be dealt with in the same way.

\bigskip 

The sheet $S'$  does not intersect $R$ nor $\theta(R)$. On the other hand, $R$ and $\theta(R)$ intersect in $\pm \OO_{(2^{n-1},1^2)}$, hence the stratum containing them is not smooth.

\bigskip

Let us deal with $R$. The Weyl group element associated with it is 
$w_R=s_{\varepsilon_1+\varepsilon_2}s_{\varepsilon_3+\varepsilon_4}\cdots s_{\varepsilon_{n-2}+\varepsilon_{n-1}}=w_0w_\Pi$ for  $\Pi=\{\alpha_1,\alpha_3,\ldots,\alpha_{n-2}\}$. In this case $\theta(w_R)=w_R$. Let us consider the injective morphism $\iota\colon SO_{2h}(k)\to SO_n(k)$ given by $\left(\begin{smallmatrix}
A&B\\
C&D\end{smallmatrix}\right)\mapsto \left(\begin{smallmatrix}
A&&B\\
&1\\
C&&D\\
&&&1\end{smallmatrix}\right)$. Then, for $\w_S$ as for $n=2h$, we choose $\w_R:=\iota(\w_S)$ and we get $V_R=\iota(V_S)$. 

Thus, for  $i=1,\,\ldots, h$ and $\epsilon_i=\pm1$,  $\w_{R} (Z(L_{\Pi})\cap T^{w_{R}})V_{R}$ is the disjoint union of the sets of matrices of the form
$x(E,D,\zeta)=\left(\begin{smallmatrix}
0&&E\\
&\zeta\\
E&& D\\
&&&\zeta^{-1}\end{smallmatrix}\right)$ with $x_i\in k$ for $i=1,\ldots, h$, $\zeta\in k^*$ and $E$, $D$ as for $n=2h$. A matrix $x(E,D,\zeta)$ lies in $S$ only if there exists $\lambda\in k^*$ such that ${\rm rk}(x(E,D)-\lambda I)\leq n$. This is possible only if $D=(\lambda+\lambda^{-1})I$ and $\zeta=\lambda^{\pm1}$. Conversely, if this is the case, a direct verification, making use of the computations for $n=2h$ and the sheet $S$, shows that $x(E, (\lambda+\lambda^{-1})I,\zeta)$ lies in $R$ if $\zeta =\lambda$ and it lies in $\theta(R)$ otherwise. Thus, set theoretically, $R\cap {\mathcal S}_{w_R}$ is a disjoint union of $2^h$ affine lines. We conclude as in the previous cases.

\subsection{Exceptional groups}

There are no non-trivial sheets of spherical conjugacy classes in types $E_8,F_4,$ and $G_2$, 
so strata of spherical conjugacy classes consists of finitely many classes, hence they are smooth. 
Let us analyze the cases for $G$ of type  $E_6$ or $E_7$. 
\begin{enumerate}
\item[$E_6$] Let $\omega\in k$ be a primitive fourth root of $1$ and let $\zeta$ be a primitive third root of $1$.  
For $a\in k^*$, let $p_{2,a}=h_1(a^2)h_2(a^3)h_3(a^4)h_4(a^6)h_5(a^5)h_6(a^4)$ and let $\OO_{2A_1}$ 
be the unipotent 
conjugacy class in $G$ of type $2A_1$. Then the only non-trivial sheet containing spherical classes is 
$$S=\left(\cup_{a\in {k},\; a^3\neq0,\,1}G\cdot p_{2,a}\right)\cup \left(\cup_{z\in Z(G)}z \OO_{2A_1}\right).$$
By \cite[Theorem 3.6]{gio-pacific}, if $a\neq b$ then  $p_{2,a}$ is not conjugate to $p_{2,b}$. 

In this case, $w_S= w_0w_\Pi$ for $\Pi=\{\alpha_3, \alpha_4,\alpha_5\}$, so $w_S=s_\beta s_\gamma$  where 
$\beta=\alpha_1+2\alpha_2+2\alpha_3+3\alpha_4+2\alpha_5+\alpha_6$ 
is the highest root and $\gamma=\alpha_1+\alpha_3+\alpha_4+\alpha_5+\alpha_6$ is the highest root in 
$\Phi\cap \beta^\perp=\Phi\cap \{\alpha_1,\alpha_3,\alpha_4,\alpha_5,\alpha_6\}$. 

We compute the set theoretical intersection $S\cap{\mathcal S}_{w_S}$ by detecting $\OO\cap{\mathcal S}_{w_S}$ 
for each orbit in $S$.  

Let us use a parametrization $x_{\pm\alpha}(\xi)$ of the root subgroups $X_{\pm\alpha}$, for $\alpha\in\{\beta,\gamma\}$ and 
$\xi\in k$, satisfying  $x_\alpha(1)x_{-\alpha}(-1)x_\alpha(1)=n_\alpha$, 
with $n_\alpha$ commuting with roots subgroups associated
with roots in $\pm\Pi$. We choose $\w_S:=n_\beta n_\gamma$.

We first consider $\OO_a=G\cdot p_{2,a}$ for $a^3\neq0,1$. Since $\beta(p_{2,a})=\gamma(p_{2,a})=a^3$ we have, for 
$x=\frac{1}{a^{-3}-1}$:
$$x_{-\beta}(x)x_{-\gamma}(x)p_{2,a}x_{-\gamma}(-x)x_{-\beta}(-x)=x_{-\beta}(-1)x_{-\gamma}(-1)p_{2,a}$$
and 
$$
\begin{array}{l}
z_a:=x_\gamma(1)x_\beta(1)x_{-\beta}(-1)x_{-\gamma}(-1)p_{2,a}x_\beta(-1)x_\gamma(-1)\\
=\w_S p_{2,a}x_\beta(-a^{-3}-1)x_\gamma(-a^{-3}-1)\in w_S TU^{w_S}\cap \OO_a.
\end{array}
$$
For $a\in k^*$ let  $b,c\in k$ satify $b^4=c^4=a^3$. 

Conjugation of $z_a$ by $h_\beta(b)h_\gamma(c)$ gives:
$$\begin{array}{l}
y_{a,b^2,c^2}:=h_\beta(b)h_\gamma(c)z_a  h_\gamma(c)^{-1} h_\beta(b)^{-1}=\\
 \w_S h_1(a^{2}(bc)^{-2})h_3(ac^{-2})h_4(c^2b^{-2}) h_5(a^2c^{-2})
 h_6(ab^2c^{-2})\cdot \\
 x_\beta(-b^2(a^{-3}+1))x_\gamma(-c^2(a^{-3}+1))\in \OO_a\cap {\mathcal S}_{w_S}.
  \end{array}
$$
which depends on $a,b^2,c^2$, for $c^2=\pm b^2$. Since $\OO_a\cap{\mathcal S}_{w_S}$ is a single 
$\Gamma_{w_S}$-orbit and  $\Gamma_{w_S}= \langle h_\beta(\omega), h_\gamma(\omega)\rangle$, we have  
$$\left(\bigcup_{a^3\neq0,1}\OO_a\right)\cap{\mathcal S}_{w_S}=\bigcup_{\epsilon=\pm1}
\left(\bigcup_{a^3\neq 0,1;\,a^3=d^2}y_{a,d,\epsilon d}\right).$$
We analyze now the orbits in $Z(G)\OO_{2A_1}$. We recall that $Z(G)=\langle p_{2,\zeta}\rangle$.
A representative of $\OO_{2A_1}$ is $u=x_{-\beta}(-1)x_{-\gamma}(-1)$, so for $0\leq l\leq 2$, the element
$$y_l:=x_{\gamma}(1)x_{\beta}(1) p_{2,\zeta^l} u x_{\beta}(-1)x_{\gamma}(-1)=\w_S p_{2,\zeta^l} x_{\beta}(-2)x_{\gamma}(-2)$$ 
lies in $p_{2,\zeta^l}\OO_{2A_1}\cap {\mathcal S}_{w_S}$.
All other elements in this set are obtained by $\Gamma_{w_S}$-conjugation:
$$
\begin{array}{l}
 y_{\zeta^{l},i,j}:=h_\beta(\omega^i) h_\gamma(\omega^j) y_l h_\gamma(\omega^{-j}) h_\beta(\omega^{-i})\\
 =\w_S h_1(\zeta^{-l}(-1)^{i+j})h_3(\zeta^{l}(-1)^{j})h_4((-1)^{i+j}) h_5((-1)^{j}\zeta^{-l})
 h_6((-1)^{i+j}\zeta^l)\cdot \\
 x_\beta(-2(-1)^{i})x_\gamma(-2(-1)^{j})
 \end{array}$$
 hence
 $$\left(\bigcup_{z\in Z(G)}z \OO_{2A_1}\right)\cap {\mathcal S}_{w_S}=
\bigcup_{\epsilon=\pm1}\left(\bigcup_{a^3=1;\,1=d^2}y_{a,d,\epsilon d}\right).$$
By Remark \ref{rem:francesco}
$S\cap{\mathcal S}_w$ is the union of two disjoint irreducible components, each isomorphic to the image of the curve 
$x^3=y^2$, for $x,y\neq0$, through the morphism $(x,y)\mapsto (x^{-1},xy^{-1},x^2y^{-1},x,y(x^{-3}+1))$.

\item[$E_7$] For $a\in k^*$, let 
$q_{3,a}=h_1(a^2)h_2(a^3)h_3(a^4)h_4(a^6)h_5(a^5)h_6(a^4)h_7(a^3)$ and let $\omega$ be a fourth primitive root of $1$. 
Let $\OO_{3A_1''}$ be the unipotent conjugacy class in $G$ of type $3A_1''$. 
Then
$$S=\left(\cup_{a\in {\mathbb C},\; a\neq0,\pm1}G\cdot q_{3,a}\right)\cup \left(\cup_{z\in Z(G)}zC_{3A_1''}\right)$$
is the only non-trivial sheet containing spherical classes. Here, $w_S=w_0w_\Pi$ for $\Pi=\{\alpha_2,\alpha_3,\alpha_4,\alpha_5\}$ so $w_s=s_\beta s_\gamma s_7$ for 
$\beta=2\alpha_1+2\alpha_2+3\alpha_3+4\alpha_4+3\alpha_5+2\alpha_6+\alpha_7$ the highest root and $\gamma=\alpha_2+\alpha_3+2\alpha_4+2\alpha_5+2\alpha_6+\alpha_7$ the highest root in $\Phi\cap\beta^\perp$. We choose $\w_S:=n_\beta n_\gamma n_{\alpha_7} $ and we will argue as we did for $E_6$.

Let us first consider $\OO_a=G\cdot q_{3,a}$ for $a^2\neq0,1$. Since $\beta(q_{3,a})=\gamma(q_{3,a})=\alpha(q_{3,a})=a^2$ we have, for 
$x=\frac{1}{a^{-2}-1}$:
$$\begin{array}{l}
x_{-\beta}(x)x_{-\gamma}(x)x_{-\alpha_7}(x)q_{3,a}x_{-\alpha_7}(-x)x_{-\gamma}(-x)x_{-\beta}(-x)\\
=x_{-\beta}(-1)x_{-\gamma}(-1)x_{-\alpha_7}(-1)q_{3,a}
\end{array}$$
and 
$$
\begin{array}{l}
z_a:=x_\gamma(1)x_\beta(1)x_{\alpha_7}(1)x_{-\beta}(-1)x_{-\gamma}(-1)x_{-\alpha_7}(-1)q_{3,a}x_{\alpha_7}(-1)x_\beta(-1)x_\gamma(-1)\\
=\w_S q_{3,a}x_\beta(-a^{-2}-1)x_\gamma(-a^{-2}-1)x_{\alpha_7}(-a^{-2}-1)\in w_S TU^{w_S}\cap \OO_a.
\end{array}
$$
For $a\in k^*$, let  $b,c,d\in k$ satify $b^4=c^4=d^4=a^2$. 

Conjugation of $z_a$ by $h_\beta(b)h_\gamma(c)h_{\alpha_7}(d)$ gives:
$$\begin{array}{l}
y_{a,b^2,c^2,d^2}:=h_\beta(b)h_\gamma(c)h_{\alpha_7}(d)z_a  h_{\alpha_7}(d)^{-1}h_\gamma(c)^{-1} h_\beta(b)^{-1}\\
= \w_S 
 h_2(ac^{-2})h_3(a^2b^{-2}c^{-2}) h_5(ab^{-2})
 h_7(a^3b^{-2}c^{-2}d^{-2})\cdot \\
 x_\beta(-b^2(a^{-2}+1))x_\gamma(-c^2(a^{-2}+1))x_{\alpha_7}(-d^2(a^{-2}+1))\in \OO_a\cap {\mathcal S}_{w_S}.
  \end{array}
$$
which depends on $a,b^2,c^2,d^2$, for $c^2=\pm b^2=\pm d^2=\pm a$.
As all  elements in  $w_S T^{w_S}U^{w_S}\cap \OO_a$ form a single orbit for the group 
$\Gamma_{w_S}= \langle h_\beta(\omega), h_\gamma(\omega),h_{\alpha_7}(\omega)\rangle$, we have
$$\left(\bigcup_{a^2\neq0,1}\OO_a\right)\cap{\mathcal S}_{w_S}=\bigcup_{\epsilon,\eta,\theta=\pm1}
\left(\bigcup_{a^2\neq 0,1}y_{a,\epsilon a,\eta a,\theta a}\right)$$ 
and, for $b^2=\epsilon a$, $c^2=\eta a$, $d^2=\theta a$ we have
$$
\begin{array}{l}y_{a,\epsilon a,\eta a,\theta a}=\w_S  h_2(\eta)h_3(\epsilon\eta) h_5(\epsilon)
 h_7(\epsilon\eta\theta)\cdot \\
 x_\beta(-\epsilon(a^{-1}+a))x_\gamma(-\eta(a^{-1}+a))x_{\alpha_7}(-\theta(a^{-1}+a)).\end{array}$$

Let us now consider the orbits in $Z(G)\OO_{3A_1''}$. We recall that $Z(G)=\langle q_{2,-1}\rangle$.
The class $\OO_{3A_1}$ is represented by $u=x_{-\beta}(-1)x_{-\gamma}(-1)x_{-\alpha_7}(-1)$, so, for $\xi=\pm1$, the element
$$
\begin{array}{l}
y_\xi:=x_{\gamma}(1)x_{\beta}(1)x_{\alpha_7}(1) q_{3,\xi} u x_{\alpha_7}(-1)x_{\beta}(-1)x_{\gamma}(-1)\\
=\w_S q_{2,\xi} x_{\beta}(-2)x_{\gamma}(-2)x_{\alpha_7}(-2)\in q_{3,\xi}\OO_{3A_1''}\cap {\mathcal S}_{w_S}.
\end{array}$$ 
All other elements in this set are obtained by $\Gamma_{w_S}$-conjugation:
$$
\begin{array}{l}
h_\beta(\omega^i) h_\gamma(\omega^j) h_7(\omega^l) y_\xi h_7(\omega^{-1})h_\gamma(\omega^{-j}) h_\beta(\omega^{-i})\\
 =\w_S h_2(\xi(-1)^{j})h_3((-1)^{i+j})h_5(\xi(-1)^{i}) h_7(\xi(-1)^{i+j+l})\cdot \\
 x_\beta(-2(-1)^{i})x_\gamma(-2(-1)^{j})x_{\alpha_7}(-2(-1)^{l})
 \end{array}.$$
 We conclude that
 $$S\cap{\mathcal S}_{w_S}=\bigcup_{\epsilon,\eta,\theta=\pm1}
\left(\bigcup_{a^2\neq 0}y_{a,\epsilon a,\eta a,\theta a}\right),$$ 
which, by Remark \ref{rem:francesco} is isomorphic to a disjoint union of $8$ copies of an affine line, with coordinate ring $k[a+a^{-1}]$.

\end{enumerate}

\section*{Acknowledgements}

The present work was partially supported by Progetto di Ateneo CPDA125818/12 of the University of Padova.

\end{document}